\documentclass[12pt,reqno]{amsart}

\usepackage{amstext,amsmath,amssymb,amsthm}
\usepackage[dvips]{graphics,color}

\begin{document}

\def\theequation{\thesection.\arabic{equation}}
\newcommand{\cal}{\mathcal}
\renewcommand{\rq}[1]{(\ref{#1})}
%\label{De
\newcommand{\ov}{\overline}
\renewcommand{\v}{ {\,\bf v }}
\newcommand{\w}{ {\,\bf w }}
\newcommand{\diam}{\mbox{diam}}
\newcommand{\dsize}{\displaystyle}
\newcommand{\supp}{\mbox{supp}}
\renewcommand{\Re}{\mbox{Re\,}}
\renewcommand{\Im}{\mbox{Im\,}}
\renewcommand{\supp}{\mbox{supp}}
\newcommand{\R}{{ \mathbb R  }}
\newcommand{\C}{ \mathbb C }
\newcommand{\N}{ \mathbb N }
\newcommand{\Z}{ \mathbb Z }

\newcommand{\BEL}{\begin{equation}\label}
\newcommand{\EE}{\end{equation}}
\renewcommand{\rq}{\eqref}
\renewcommand{\L}{\mathscr{L}}

\renewcommand{\t}{\tilde}
\renewcommand{\medskip}{\vskip .5 cm}
\newtheorem{Thm}{Theorem}[section]
\newtheorem{Lemma}[Thm]{Lemma}
\newtheorem{Cor}[Thm]{Corollary}
\newtheorem{Prop}[Thm]{Proposition}

\begin{center}
{\Large \bf Split functions, Fourier transforms and multipliers}
\medskip
Laura De Carli and Steve Hudson
\end{center}
\bigskip
\centerline{\bf Abstract}

 \medskip We study the effect of a splitting operator $S_t$ on the  $L^p$ norm of the  Fourier transform of a function $f$ and on the operator  norm of a  Fourier multiplier $m$.   Most of our results assume $p$ is an even integer, and are often stronger when $f$ or $m$ has compact support.

\medskip\noindent {\bf Mathematics Subject Classification:} 42B15, 42B10

\section*{1. Introduction}
\setcounter{section}{1} \setcounter{Thm}{0} \setcounter{equation}{0}

We study the effect of a splitting operator $S_t$ on the Fourier transform of a function $f$.
For $f\in L^2(\R,\C)\cap L^1(\R,\C)$, let  $\tau_t f(x)= f(x -t)$. Let
$f_+(x)= f(x)\chi_{\{x > 0\}}$ and $f_-(x)= f(x)\chi_{\{x < 0\}} $,
where $\chi$ denotes the characteristic function; so $f=f_++f_-$ almost everywhere.
For $t>0$, let
$$S_t f (x)= \tau_tf_+ + \tau_{-t}f_-$$
We call $S_t$ a {\it splitting operator} on $L^2 \cap L^1(\R,\C)$.

For functions of  $x=(x', x_n) \in \R^{n-1} \times \R = \R^n$, we will often treat $x'$ as a parameter and
view $f$ a function of the single variable $x_n$. Most of our translations and convolutions will be in this variable. For $f\in L^2(\R^n,\C)\cap L^1(\R^n,\C)$, let $\tau_t f(x)= f(x', x_n -t)$,  $f_+(x)= f(x)\chi_{\{x_n\ge 0\}}$,  $f_-(x)= f(x)\chi_{\{x_n\leq 0\}} $ and  $S_t f (x)= \tau_{t}f_+  + \tau_{-t}f_- $. The Fourier transform of $S_t f$  is
\begin{equation}\label{FofSplit}
 \widehat{S_tf}(y)= \int_{\R^n } e^{-2\pi i xy} S_t f(x)dx= e^{-2\pi i t y_n} \hat f_{+}(y)  + e^{2\pi i t y_n} \hat f_{-}(y).\end{equation}

We study how $N_tf = ||\widehat{S_tf}||_p = (\int_{\R} |\widehat{S_tf}(x)|^p\ dx )^{1/p}$ depends on $t$.
The results are often trivial when $p=2$, since by the  Plancherel theorem, $N_tf$ is independent of $t$.
Throughout the paper, $1<p<\infty$ is fixed. In Sections 2 and 3, we will also assume that $p\ge 2$ and is an even integer.
Our main result is for real-valued functions in the class
\begin{equation}\label{S}
S=\{ f\in L^2\cap L^1(\R) \ : \
\quad f_+*f_-(x)\ \ \mbox{decreases for} \  x\ge 0\}
\end{equation}
where by $L^p(\R^n)$ we mean $L^p(\R^n,\,\R)$.
We will show in Section 3 that many standard functions satisfy \eqref{S}.
For example, non-negative even functions that decrease for $x>0$, (also called {\it radially-decreasing functions} or simply {\it bump functions}) satisfy \eqref{S}. Our main result is the following.

\begin{Thm}\label{Decr}
Assume $f\in S$ is even and non-negative, and that $p$ is even.
Then, $N_tf$ decreases with $t$.
In particular $N_t f\leq N_0 f$.
\end{Thm}

 We prove this in Section 3 and show that the hypothesis $f\in S$ cannot be removed.
 It is unclear whether the other assumptions on $f$ are necessary. It is also unclear whether this theorem generalizes to functions in $\R^n$, except in some simple cases (See Corollary \ref{rect}). In Section 2, we prove that $N_t f$ is eventually constant when $f$ has compact support.  \begin{Thm}\label{Const}   Suppose $h=f+ig \in L^2 (\R, \C)$ is supported on $[-A,A]$
and $p$ is even.  Then,  $N_t h=N_{t_0}h$ whenever
\begin{equation}\label{lower} t\ge  t_0 = \frac{(p-2)A}{4}.\end{equation}
\end{Thm}

The lower bound in \rq{lower} cannot be improved in any obvious way.
For example, let $f(x)=\chi_{(-1,1)}(x)$ and $p=4$. Then
\BEL{FTS} \widehat {S_t f}(y) = 2\int_t^{1+t} \cos(2\pi x y)dx
   = \frac{   \sin (2\pi(1+t) y)-\sin (2\pi t y) }{ \pi y}\EE
and
$$(N_t f)^4 = ||S_tf||_4^4=\frac{1}{24}\left( 6+ (1-2t)^3 + |1-2 t|^3) \right)$$
which is constant for $t\ge t_0 = \frac 12$, but not for $t< \frac 12$.

We offer numerical evidence (but no rigorous proof) that
Theorem \ref{Const} can fail when $p=3$. Let $f(x)= \chi_{(-1,1)}(x)$.
If Theorem \ref{Const} were true for $p=3$, then $N_t f$ would be constant for $t\ge \frac 14 $. Numerical integration gives $(N_{\frac 14}f)^3 \sim  2.6247$,  $(N_{1}f)^3 \sim   2.6124$, $(N_{5}f)^3 \sim 2.6116$ and $(N_{12}f)^3 \sim 2.6121$. If this data is accurate, it also shows that $N_tf$ is not decreasing, so that Theorem \ref{Decr} also fails when $p=3$.

\medskip
In Section 4 we study the effect of the splitting operator $S_t$ on the norm of a Fourier multiplier.
Let $m\in L^\infty (\R^n, \C)$. Suppose the operator  $T_mf(x)= \int_{\R^n} e^{2\pi i xy} m(y) \hat f(y) dy$,
initially defined for  $f \in C^\infty_0(\R^n)$, extends to a bounded operator on $L^p(\R^n )$.
Then we say that $m$ is a $(p,p)$ Fourier multiplier
and we write its operator norm (or {\it multiplier norm}) as $|||T_m|||_{p,p}^\R $ or simply by $|||m|||_{p,p}^\R $.
If we apply $T_m$ to complex-valued functions, it has a possibly larger operator norm, denoted $|||m|||_{p,p}$. Indeed, $|||m|||_{p,p}^\R\leq |||m|||_{p,p}  \leq 2|||m|||_{p,p}^\R$. We relate $|||S_t m|||_{p,p}$ to   $|||m_+|||_{p,p}$ and the norm of the half line multiplier  $s = \chi_{(0, \, \infty)}$ (see  Theorems \ref{Mult}, \ref{HalfLine} and Proposition \ref{T-Poly}). We also estimate the $(p,p)$ norm of $m_+$, with exact values in some  cases (See Corollary
\ref{P-extram+}).
\medskip

One of the authors became interested in splitting  operators while  investigating $L^p$ Lebesgue constants,  i.e.
the functions  $L(N, W,p)=\int_{[0,1]^{n}}\left\vert
   \sum_{j \in \Z^n\cap N W}e^{2\pi i j\cdot x}  \right\vert ^{p}dx$, where  $W$ is a convex set in $\R^n$,  $p\ge 1$ and   $N>0$.   This problem has a long history (see \cite{L}, \cite{NP}, \cite{Y}; see also \cite{TB} and the references cited there). These constants depend on $N$ in an explicit way,   but  their dependence on  $W$ is not as clear. It would be interesting to  estimate  the  $L^p$ Lebesgue constants  when  $W=P_{k}$ is a regular polygon in $\R^2$ with $2^k$ sides,
   especially with values of $p$ for which these  constants
   are not uniformly bounded in terms of $k$  (see \cite{AD}).
   To do this, it might be useful to compare the $L^p$ norms  of $S_t( P_k)$ and that of $P_k$  when   $S_t$  divides   $\chi_{P_{k}}$ into  two congruent parts.

The comparison of the   $(p,p)$ norms of $\chi_{P_{ k}} $   and  $S_t(\chi_{P_{ k}})$ is also very interesting. When   $\frac 43<p<4$, asymptotic estimates of the    $(p,p)$ norm of $\chi_{P_{ k}}$  in terms of $k$ are in  \cite{Co1}, (see also \cite{Co2}, \cite{CS}, and see \cite{F} for the multiplier problem for the disk). But we know of no such estimates for other $p$, or when $P_k$  is not convex - for example, when $P_k$ is the disjoint union of two congruent polygons. See Theorem \ref{Mult} and Proposition \ref{T-Poly} for results in this direction.

 \medskip
{\it Acknowledgement:} The authors wish to thank Marshall Ash for pointing out related problems which inspired this work.
\medskip

\section*{2. Proof of Theorem \ref{Const}   and related results}
\setcounter{section}{2} \setcounter{Thm}{0} \setcounter{equation}{0}

Let $h$ be as in Theorem \ref{Const}. We plan to show that  $(N_t h)^p = C+ \sum H_j(t)$, where $C$ is a constant and each $H_j$ is a convolution of $p$ functions with compact support (see formula \rq{N}). This implies, using Lemma \ref{Convo} below, that each $H_j$ has compact support, and that $N_th = C$ for large $t$, as desired.

\begin{Lemma}\label{Convo}
Let $\phi_1$  and $\phi_2$ be integrable functions supported  in $[a,b]$ and $[c,d]$, resp. Then
$\phi_1*\phi_2(x)= \int_a^b \phi_1(y) \phi_2(x-y) dy$ is supported in
$[a+c,b+d]$.
\end{Lemma}
\medskip
\noindent {\it Proof.}
Suppose $x> b+d$.
If $y\le b$, then $x-y > d$ and $\phi_2(x-y)=0$, so $\phi_1*\phi_2(x)=0$.
The case $x< a+c$ is similar.
$\Box$
\medskip
\noindent
{\it  Proof of Theorem \ref{Const}.} Let $h= f+ig \in L^2(\R,\C)$ with  support in $[-A,A]$ and let $t\ge t_0$.
Denote the even part of $f$ by $f_e(x)=\frac 12(f(x)+ f(-x))$ and the odd part by $f_o = f-f_e$. Note that if $\phi$ is real, then $\Re \hat \phi = \widehat {\phi_e}$. With
$\phi = S_tf$, we get $\Re  \widehat {S_tf} (y)   =  \widehat {(S_tf)_e} (y) = \widehat {S_t (f_e)} (y)$.
Likewise, $\Im \hat \phi =  -i \hat \phi_o$ and $ \Re \widehat {S_t (ig)} (y) = -\Im \widehat {S_t g} (y) =
i \widehat {S_t (g_o)} (y)$. By linearity,
$$
S_t h (x)=    S_t(f _e) (x) + S_t(f _o) (x) + i(S_t (g_e)(x)+ S_t (g_o) (x))
 .
$$
So,
$$ \Re \widehat {S_t h}(y)=  \widehat {S_t(f_{e})}(y) +i  \widehat {S_t (g_{o})}(y). $$
By \rq{FofSplit},
\begin{equation}\label{ReS}
\Re \widehat {S_th}  = e^{-2\pi i t y}\left(\hat f_{e+} +i\hat g_{o+}\right)+ e^{2\pi i t y}\left(\hat f_{e-} +i\hat g_{o-}\right)
. \end{equation}
Im $\widehat {S_t h}$ can be written as a similar sum, but with terms containing   $f_{o\pm}$ and $g_{e\pm}$.
Let $h_1$ through $h_4$ be the functions $f_{e-}$, $f_{o-}$, $ig_{e-}$, $ig_{o-}$, which are supported on $[-A,0]$.
Let $h_5$ through $h_8$ be the functions  $f_{e+}$, $f_{o+}$, $ig_{e+}$, $ig_{o+}$, which are supported on $[0,A]$.
Squaring both sides of \rq{ReS}  and the similar formula for Im $\widehat {S_th}$, and adding them,
\begin{equation}\label{sumprod} |\widehat {S_th }|^2 = \sum e^{2\pi i t y C_\alpha}\prod_{j=1}^8 \hat h_j^{a_j}
\end{equation}
where the sum is over certain $\alpha = (a_1, \ldots a_8)$, such that each $a_i$ is a non-negative integer,
and $\sum a_i = 2$. We do not assert that all such $\alpha$ appear exactly once in this sum.
Here, $C_\alpha =\sum_{j=1}^4 a_j - \sum_{k=5}^8 a_k$.
Raising both sides of \eqref{sumprod} to the power $\frac p2$, we get a sum of the form
$$|\widehat {S_th}|^p = \sum_\beta k_\beta e^{2\pi i t y P_\beta}\prod_{j=1}^8 \hat h_j^{b_j}$$
where $\beta = (b_1, \ldots b_8)$, each $b_j$ is a non-negative integer, $\sum b_j = p$,
and $ P_\beta=\sum_{j=1}^4 b_j - \sum_{k=5}^8 b_k$, which we will write as $p_+-p_-$.

When we integrate this, the exponential factors produce inverse Fourier transforms,
which convert products to convolutions. So,
\BEL{N}
(N_t  h) ^p = \int_\R |\widehat {S_th}(y)|^p dy= \sum_\beta k_\beta H_\beta (P_\beta t)
\EE
where $H_\beta$ is a convolution of $p$ functions (the ones appearing in ${\dsize \prod_{j=1}^8 }{\hat h_j}^{b_j}$).
Since $  P_\beta  =p-2p_-$,  it is even. Whenever $P_\beta=0$, the summand in \rq{N} is constant in $t$.
If $P_\beta\ge 2$, then $2p_- \le p-2$ and $P_\beta  t \ge P_\beta  t_0 \ge A(p-2)/2 \ge Ap_-$.

By Lemma \ref{Convo}, $H_\beta$ is supported on $[-Ap+,\ Ap_-]$, so $H_\beta(P_\beta t) =0$.
For summands with $P_\beta\le -2$,
we have  $p-2 \ge 2p_+$, leading to the same conclusion.
This proves that $(N_t  h) ^p  $ is constant for $t\ge t_0$. $\Box$

\medskip
 Theorem \ref{Const} has an analogue for Fourier series, with a similar proof, which we omit.
 Let  $f = \{c_k\}_{k= -A}^A$ be a finite sequence of complex numbers.   Let $\hat f(x) =\sum_{k = -A}^{+A} c_k e^{-2\pi ikx}$ be the corresponding trigonometric polynomial.
Let  $p$ and $t$ be positive integers,
with $p$ even. Define $S_tf = \{b_k\}_{|k| \le A+t}$, where $b_{k-t} = c_{k}$ for $k<0$, $b_{k+t} = c_{k}$ for $k>0$,
$b_{-t} = b_{t} = c_{0}/2$ and $b_{k}=0$ otherwise. Note that
\BEL{shiftFT}\widehat{S_tf}(x) = c_0 \cos(2\pi tx) + \sum_{k = 1}^{+A} c_k e^{-2\pi i(k+t)x}+ \sum_{k = -A}^{-1} c_k e^{-2\pi i(k-t)x}.   \EE
\medskip
\begin{Thm}\label{ConstSeries}   For $t\ge t_0 = \frac {(p-2)A}{4}$, $N_tf = ||\widehat{S_tf}||_{L^p[0,1]}$ is independent of  $t$.
\end{Thm}

The split $f= f_- + f_+$ defined in the Introduction is a bit arbitrary and can be generalized.
Suppose $f_1$ is supported in $[-A, b]$ and $h_2$ is supported in $[ -b, A]$, with $|b| \leq A$ (so, the supports overlap if $b>0$).
Let $f= f_1+ f_2$ and let $\mathcal S_tf(x)= f_2(x-t)+ f_1(x+t)$ for $t\in \R$.
Let $N_t f= ||\widehat{S_tf} ||_p$.

\begin{Cor}\label{GenConst}  With $f$, $\mathcal S_tf$ and $N_t$ as above, and $p$ even, let $t \ge  t_0 = \frac{p-2}{4}(A+b) +b$.
Then  $N_tf=N_{t_0}f$.
\end{Cor}
\noindent
{\it Proof.} Let $t > \frac{p-2}{4}(A+ b)+  b$.
  Note that $g= \mathcal S_{b}f$ is supported in $[-A-b, A+b]$ and $S_t g= \mathcal S_{t+b} f$.  By Theorem \ref{Const}, $N_{t} g$ is constant if $t > \frac{p-2}{4} (A+b)$. So, $ N_{t }f= N_{t-b }g$ is constant if $ t > \frac{p-2}{4}(A+ b)+  b   $ as required. $\Box$

\medskip
One can generalize this result  a bit more, and similarly Theorem \ref{ConstSeries}, with essentially the same proof.
Suppose $f_1$ is supported in $[-A, b_1]$ and $f_2$ is supported in $[b_2, A]$, with $|b_1|, |b_2|  \leq A$.
If $ t\! \ge \! t_0 = \frac{p-2}{4}\!\left(A+ \frac{|b_1+b_2|}{2}\right)$ $ +\frac{p-2}{8} (b_1-b_2)$, then $N_t f=N_{t_0}f$.

\medskip
Theorem \ref{Const} also has an n-dimensional analogue.
\medskip
\begin{Cor}\label{Constn}  Assume  $ f(x)\in L^1 \cap L^2(\R^n, \C)$ is
supported where $|x_n|<A$. If $p$ is even, and $ t\ge  t_0 = \frac {(p-2)A}{4}$, then $N_tf$ is independent of $t$.
\end{Cor}
\noindent
{\it Proof.} Note that $\widehat{f}(y) = {\cal F}_{n-1} {\cal F}_{1} f  (y)$, where
${\cal F}_{n-1}$ denotes the Fourier transform of $f(x', x_n)$ with respect to $x'$, with $x_n$ fixed,
and ${\cal F}_{1}$ denotes the Fourier transform with respect to $x_n$, with $x'$ fixed.
 We will sometimes write $\hat g$ as ${\cal F}_{1} g$ for functions $g$ defined on $\R$.
Fix $y'\in R^{n-1}$ and $t\ge t_0$. Define $g(x_n)= g(y',x_n)
=  {\cal F}_{n-1}f(y', x_n)$, and $g_+$, $g_-$, $S_tg $ as usual.
Note that  $S_t$ commutes with ${\cal F}_{n-1}$ because they act on different variables, so
${\cal F}_{1} { S_t g}(y_n) = {\cal F}_{1}{\cal F}_{n-1}S_t f(y',y_n)= \widehat {S_tf}(y).$
So,
\begin{equation}\label{intg} \int_{\R} |\widehat {S_tf} (y)|^p \ dy_n =
\int_{\R}| {\cal F}_{1}(S_tg)(y_n)|^p dy_n = (N_t g)^p ,\end{equation} which is constant for $t\ge t_0$, by
Theorem \ref{Const}.
Since $S_t$ is an isometry on $L^1$ and $L^2$,  the Hausdorff-Young inequality implies that $\widehat{S_t f} \in L^2\cap L^\infty \subset L^p(\R^n)$. Integrating \eqref{intg} over $y'$ concludes the proof.
$\Box$
\vskip 10pt

In the corollary, we split $\R^n$ around the hyperplane $H= \{x: x_n=0\}$ only for convenience. Since all the norms are invariant under rotations and translations, we could use any other hyperplane in $\R^n$ instead. After adjusting the definition of $S_tf$ in an obvious manner, an analogue of Corollary \ref{Constn} holds for any such $H$.

\medskip

The next result will be applied to the study of multipliers in Section 4, and involves functions that depend on $x$ and $t$. Let $F(x,t):\R^n\times [0,\ \infty)\to\C$  such that     $F(\cdot,t)\in L^1 \cap L^2(\R^n,\C)$   for every $t\ge 0$. We let   $F_+(x,t)= F_+(x', x_n, t)= F(x,t)\chi_{ \{x_n>0\} }(x)$ and similarly for $F_-$.
We let $S_t F(x,t)=   S_t F(x',\, .,\, t)= F_+(x', x_n-t,t) + F_-(x', x_n+t,t)$. That is, $S_t$ acts only on the  $x_n$ variable.

\begin{Lemma}\label{CWt}
Suppose that $F(x,t)$ as above is supported where $|x_n|<A$  and, for every  $t\ge 0$, that $\widehat {S_tF}(y, t)$ is real-valued. If $t\ge t_0=\frac{(p-2)A}{4}$ then 
\begin{equation}\label{Newt}
\int_{\R^n}|\widehat {S_tF} (y,t)|^pdy= {p\choose \frac p2}\!\! \int_{\R^n} \left(  \widehat { F}_+ (y, t)  \right)^{\frac p 2} \left(
\widehat  { F}_-  (y,t)\right)^{\frac p 2}dy.
\end{equation}
\end{Lemma}

\medskip
\noindent
{\it Proof.} Fix $y'\in\R^{n-1}$ and $t>t_0$. Let $f(x_n)= {\cal F}_{n-1} F(y', x_n, t)$.  Observe  that   $S_t$ commutes with ${\cal F}_{n-1}$, and so $ {S_tf}(x_n)=  {\cal F}_{n-1}  S_t  F(y', x_n, t)$. Thus, $ {\cal F}_1{S_tf}(y_n)=  \widehat{ S_t  F}(y', y_n, t) $. Similarly, ${\cal F}_1 f_{\pm}(y_n)   = \widehat  F_{\pm}(y', y_n, t)$.
We argue as in   Theorem \ref{Const}. By \rq{FofSplit} and the binomial theorem,
\begin{align}
\nonumber  &\int_\R  |\widehat{ S_t  F}(y', y_n, t)|^p  dy_n   \nonumber  =\int_\R|{\cal F}_1 {S_tf}(y_n)|^pdy_n \\\nonumber &=\int_\R (e^{-2\pi i t y_n} {\cal F}_1 f_+(y_n)+e^{2\pi i t y_n} {\cal F}_1 f_-(y_n))^p dy_n
 \\ &=  \sum_{k=0}^p {p\choose k} \int_{\R}
e^{ 2\pi i ( p-2k) t y_n} \left({\cal F}_1  f_+  (y_n)\right)^k   \left({\cal F}_1f_-(y_n)\right)^{p-k}\!dy_n \label{e-1}
\end{align}
which is an analog of \rq{N}. Reasoning as in the proof of Theorem \ref{Const}, the    integral on the right hand side of \eqref{e-1} simplifies to $(f_+)^{*k}* \left( f_-\right)^{*(p-k)} ( ( p-2k) t)=0$, unless $k=\frac p2$. The only nonzero term in \rq{e-1}
is ${p\choose \frac p2} \int_\R \left({\cal F}_1 f_+ (y_n)  \right)^{\frac p2}\left({\cal F}_1 f_- (y_n)  \right)^{\frac p2}dy_n$,  and
so
\begin{align*}
&\int_\R  |{\cal F}_1(S_t  F)(y', y_n, t)|^pdy_n  = {p\choose \frac p2} \int_{\R }({\cal F}_1 f_+( y_n )) ^{\frac p2} ({\cal F}_1 f_-( y_n )) ^{\frac p2}dy_n
\\
&={p\choose \frac p2} \int_{\R }(\widehat F_+( y',y_n,t )) ^{\frac p2} (\widehat F_-(y', y_n,t )) ^{\frac p2}dy_n
.\end{align*}
Integrating this with respect to $y'$,
\eqref{Newt} follows.
$\Box$

\bigskip
\section*{ 3. On the decreasing function $N_t$}
\setcounter{section}{3} \setcounter{Thm}{0} \setcounter{equation}{0} \setcounter{subsection}{0}

  We  prove  Theorem \ref{Decr}   and examine the class $S$ in some detail  to show that it is rather large.
  We show by example that the assumption $f\in S$ cannot be removed from the theorem.

  \subsection{Proof of Theorem \ref{Decr}}

In this section we will assume that $f$ is even and $f\in S$ unless specified otherwise. We need a few lemmas about convolution.

\medskip
\begin{Lemma}\label{ConvoEven}
 $\tau_t f_{ +}*\tau_{-t}f_{ -}(x)$ is even.
\end{Lemma}
\noindent
{\it Proof.}   Since $f$  is even, $\tau_{t}f_{+}(-x)= \tau_{-t}f_{ -}(x)$, and  we have
 \begin{align*} \tau_{t}f_{+}*\tau_{-t}f_{ -}(-x)&=
 \int_{\R}  \tau_{t}f_{+}(y) \tau_{-t}f_{ -}(-x-y) dy  \\&=    \int_{\R}  \tau_{-t}f_{ -}(-y) \tau_{t}f_{+}(x+y) dy \\  &=  \int_\R \tau_{-t}f_{ -}(y')\tau_t f_+(x-y') dy' =\tau_{t}f_{ +}*\tau_{-t}f_{ -}(x).\    \Box \end{align*}
\medskip

\begin{Lemma}
\label{adjoint}
Let $g$, $h$ and $F\in L^2(\R)$. Let $U= h*F$ and $V= g*F$.
Suppose  $F = f_1*f_2*\ldots *f_n$, where  every $f_i$ is either $f_+$ or $f_-$. Then
\begin{equation}\label{f3}
U*\t V =   h*\t g* ( f_+*f_-)^{*n}  \end{equation}
where   $\t g(x)= g(-x)$.
\end{Lemma}

\noindent
{\it Proof.} Note that $\t f_+(x)= f_-(x)$ and
$\t V = \t g* \t F$. Similarly, each $f_i* \t f_i= f_+*f_-$, so
$F*\t F = (f_+*f_-)^{*n}$ and \eqref{f3} follows. $\Box$

\medskip
\begin{Lemma}\label{ConvoDecr}   Assume that  $u \ge 0$ is supported on $(-\infty, 0]$ and that  $v$ is decreasing  on $[0,\infty)$.
Then $\dsize u*v $ is decreasing on  on $[0,+\infty)$.
\end{Lemma}
\noindent
{\it Proof.} If $0\leq x_1<x_2$, then
$$u*v(x_2)-u*v(x_1)= \int_{-\infty}^0   u(y)\left(v(x_2-y) - v(x_1-y)\right)\, dy \leq 0.\ \Box$$
\medskip
{\it Proof of Theorem \ref{Decr}:}  Assume that $f\in S$ and  $p=2m$, with $m$ a positive integer.
By the Plancherel  theorem,
\begin{eqnarray}\nonumber ( N_t f)^{2m} &=& || \widehat {S_tf}||_{2m}^{2m} = || (\widehat {S_tf})^m||_2^2= ||(S_t f)^{*m}||_2^2 \\ \label{DecrEts}&=&  || (\tau_t f_{ +} +\tau_{-t} f_{-} )^{*m}||_2^2.\end{eqnarray}
By a variation of the binomial theorem, $(\tau_t f_{ +} +\tau_{-t} f_{-} )^{*m}(x)$ is a linear combination (with positive coefficients) of terms of the form
$$(\tau_t f_{ +})^{*i}*(\tau_{-t} f_{-})^{*(m-i)}(x)= f_{+}^{*i}*f_{-}^{*(m-i)}(x-(2i-m)t) $$
where $0\le i\le m$.  Let $G_i = f_{+}^{*i}*f_{-}^{*(m-i)}$.
So, $(N_t  f)^p = N_t^{2m}(f)$ is a weighted sum of integrals of the following form:
\begin{eqnarray*} N_{ijt} &=& \int_\R G_i(x-(2i-m)t)G_j(x-(2j-m)t)dx \\ &=&  \int_\R G_i(x-2(i-j)t)G_j(x)dx= G_j* \widetilde G_i(2(i-j)t). \end{eqnarray*}
We will show that each $N_{ijt}$ decreases with $t>0$, and may assume $i\ge j$.
Let
$F(x)= f_{+}^{*j}*f_{-}^{*(m-i)}$  and $g = f_{+}^{*(i-j)}$, which is supported on $[0,\infty)$. So, $g*F(x) = G_i(x)$.
Likewise, $G_j(x) = h*F(x)$ where $h = f_{-}^{*(i-j)}$ is supported on $(-\infty,0]$. Applying Lemma \ref{adjoint} with $n=m+j-i$,
$$G_j*\t G_i = h*\widetilde  g*  [f_{+}*f_{-}]^{*(m+j-i)} = h*\widetilde  g* H . $$

By \eqref{S} and Lemma \ref{ConvoEven}, $f_{+}*f_{-} (x)$ is radially decreasing, so $H(x)$ is too.
By Lemma \ref{ConvoDecr} (with $u=h$ and $v=H$ ), $h*H(x)$ decreases for $x\ge 0$.
Applying the lemma again (with $u=\t g$), we see that
$N_{ijt} = \t g* (h*H)(2(i-j)t)$ decreases for $t\ge 0$. $\Box$.
\medskip
\noindent
{\it Example.} Theorem \ref{Decr} can fail without the assumption that $f\in S$.
Let $p=4$. The function $f(x)= \chi_{\{|x|<1\}} +  \chi_{\{10<|x|<11\}}$  is even
and non-negative, but $f_+* f_-(x)$ is not decreasing, so $f\not\in S$. The Fourier transform of $S_tf(x)=\chi_{\{ t<|x|<1+t\}}(x) +  \chi_{\{10+t<|x|<11+t\}}(x)$  can be explicitly evaluated. By the residue theorem (or Mathematica software),  when  $4<t<5$ we have $(N_t f)^4 =\frac{8}{3} \left(4 t^3-48 t^2+192 t-247\right)$, which is increasing. If $p=6$, the same $f$ leads to similar conclusions.

\medskip
The  following corollary is a variation of Theorem \ref{Decr} for dimension $n>1$.

\begin{Cor}\label{rect} Let $f(x) = g(x')h(x_n)$, where $g\in L^1\cap L^2(\R^{n-1})$ and $h\in S$.
 If $p$ is even, $N_tf $  decreases  with  $t$.
 \end{Cor}
 \medskip
 \noindent
 {\it Proof.} Since $S_t f (x) = g(x') S_t h(x_n)$, we have
 $$N_tf= ||\widehat{S_t f}||_p =  ||{\cal F}_{n-1}g ||_{L^p(\R^{n-1})} ||{\cal F}_1 (S_t h)||_{L^p(\R )}.$$
 By Theorem \ref{Decr},  $||{\cal F}_1 ( S_th)||_{L^p(\R )}$ decreases with $t$. $\Box$

\subsection{On the class $S$ }

Proposition \ref{Decr2} below shows that $S$ contains the radially decreasing integrable functions, sometimes called  bump functions,
as well as even functions with two bumps. Roughly speaking, Lemma \ref{Closed} shows that $S$ is closed in $L^\infty$. Also, $S$ is closed under the shift operator $S_t$, for $t>0$, but we leave the proof to the reader. It  is not closed under addition. For example, let $h(x)= \chi_{\{ |x|<1\}}$ and let $g(x)=\chi_{\{10<|x|<11\}}$. Proposition \ref{Decr2} shows that $h,g\in S$, but the previous example shows that $f=h+g \not\in S$.
\medskip
\begin{Prop}\label{Decr2}   Let $0\le r <\infty$.  Suppose that $f\in L^1\cap L^\infty(\R)$ is even and non-negative.
Suppose that $f_+$ is increasing on $(0,r]$, and decreasing on $[r,\infty)$.
Then $f\in S$.
\end{Prop}
If $r=0$, for example, our hypotheses reduce to the single assumption that $f_+$ is decreasing on $[0,\infty)$. Any function $f$ as in Proposition
\ref{Decr2} can be uniformly approximated by an increasing sequence of step functions, $\{g_j\}$, such as appear in Lemma \ref{Step}, with
every $||g_j||_1 \le ||f||_1$. Then Lemma \ref{Closed} proves Proposition \ref{Decr2}.
\medskip
\begin{Lemma}\label{Closed}  Suppose $\{g_j\}_{j=1}^\infty \in S\cap L^\infty(\R)$ and $\{g_j\}$ converges
uniformly to some $f \in L^1(\R)$. Suppose that every $||g_j||_1 \le C$. Then $f\in S$.
\end{Lemma}
\medskip
\noindent
{\it Proof.} Since $f$ must be bounded, it is also in  $L^2(\R)$.
Let $0\le x_1 < x_2$. It is enough to show that
\BEL{ETSS}f_+*f_-(x_2) \le f_+*f_-(x_1)+  2\epsilon (C+||f||_1)\EE
for every $\epsilon>0$. Fix $\epsilon$ and choose $j$ such that  $||g_j-f||_\infty < \epsilon$.
Since $g_j\in S$,  $g_{j+}*g_{j-}(x_2) \le g_{j+}*g_{j-}(x_1)$. In general, $||v*u||_\infty \le ||v||_1||u||_\infty$, so
\begin{align*} ||g_{j+}*g_{j-} - f_{+}*f_{-}||_\infty &\le ||g_{j+}*[g_{j-} -f_-]||_\infty + ||[g_{j+} - f_+]*f_-||_\infty\\ & \le \epsilon (C+||f||_1).\end{align*}
This proves \rq{ETSS}.\ $\Box$
\medskip
\begin{Lemma}\label{Chi}  Suppose $I_1 = [b,c]\subset [a,d] = I_2$ and $f_1(x) = \chi_{I_1}(x)$ and $f_2(x) = \chi_{I_2}(-x)$. Then $f_1 * f_2 (x)$ decreases on  $[0,\infty)$. The same result holds when $I_2\subset I_1$ instead.
\end{Lemma}
\medskip
{\it Proof.} $f_1 * f_2 (x) = \int_\R f_1(y) f_2(x-y)\ dy = \int_b^c \chi_{I_2} (y-x)\ dy $,
which is the measure of $[b,c]\cap x+[a,d]$. This is $c-b$ for $0\le x \le b-a$, then linear (and decreasing)
in $x$ until $x=c-a$, when it becomes zero. $\Box$
\medskip
\begin{Lemma}\label{Step}   Suppose $A_k$ are nested intervals with $r\in A_1 \subset A_2 \subset \ldots A_k  \subset [0,\infty)$.
Suppose $g_+(x) = \sum_{j=1}^k c_j \chi_{A_j}(x)$, where all $c_j>0$, and $g_-(-x) = g_+(x)$. Then $g_+ * g_-$ decreases on $[0,\infty)$.
\end{Lemma}
\medskip
{\it Proof.} The convolution splits into $k^2$ pairs, each of which decreases by Lemma \ref{Chi}. $\Box$
\medskip

\medskip
\section*{ 4. Split Fourier multipliers }
\setcounter{section}{4} \setcounter{Thm}{0} \setcounter{equation}{0}

\medskip
This section studies effects of the splitting operator $S_t$  on  the  norm of Fourier multipliers
(see the Introduction for definitions and notation; for the basic properties of multipliers, see \cite{S}). When $p=2$, Plancherel's theorem  implies that  every bounded function  is a  Fourier multiplier  on $L^2(\R^n,\,\C)$ and  $ |||m|||_{2,2} =||m||_\infty$.  Since $||S_t m||_\infty = ||m||_\infty$, we see that $|||S_t m|||_{2,2} = |||m|||_{2,2} $.

In general,  $||| m|||_{p,p} $ may be larger than  $||| m|||_{p,p}^\R$.
We will write $T_m \in \mathcal R$ to mean $T_m$ maps real-valued functions into real-valued functions.
This occurs, for example, when  $m$ is even and real-valued. To see this, let $f=f_e+f_o$, where $f_e$ and $f_o$ are the even and odd components. Then $T_mf= T_m f_e + T_m f_o$  is real. For such $m$, we have $||| m|||_{p,p}^\R=||| m|||_{p,p}$ (see \cite{MZ}).

Explicit formulas for multiplier norms are known only in very few cases.
The characteristic function $\sigma$
of any segment $[a,b] \subset \R^1$  (or the {\it segment multiplier}) has the same multiplier norm  as the Hilbert transform. It is:
\begin{equation}\label{segment} n_p = ||| \sigma|||_{p,p}  =
||| \sigma|||_{p,p}^\R  =\max\left\{\!\tan (\frac p{2p}), \, \cot (\frac p{2p})\!\right\}\end{equation} (see \cite{DL1}) and \cite{P}).
The $(p,p)$ norm  of the characteristic function  of the  half line    $s(x)=\chi_{(0, \infty)}(x)$ is
$c_p^\R\! = |||s|||_{p,p}^\R =\frac 12\max\left\{ \sec(\frac p{2p}) ,\, \csc(\frac p{2p}) \right\}$ (see \cite{E} and also \cite{V}). Also,
$c_p=|||s|||_{p,p}= \frac 12 \sec(\frac p{2p}) \csc(\frac p{2p}) =\csc(\frac\pi p)$ (see \cite{HV}).
Note that $c_p^\R< c_p $ and that $c_2^\R =\frac{1}{\sqrt 2}$.

Though the  operator norm of a Fourier multiplier on $L^p(\R,\,\C)$ is translation invariant,
we note that $|||\tau_{-t} s|||_{2,2}^\R\not = |||s|||_{2,2}^\R = \frac{1}{\sqrt 2}$ for $t>0$.
To see this, let $I = (-t, t)$ and let $\hat f = \chi_I$. Note that $f$ is real-valued and
$\widehat{T_{\tau_{-t} s} f}  =  \hat f$. By Plancherel, $||T_{\tau_{-t} s} f||_2= ||f||_2$, which shows that
$|||\tau_{-t} s|||_{2,2}^\R \ge 1$.

Even  when the norm of a multiplier  $m$ is known, the norm of the split multiplier $S_t m$ cannot usually be explicitly evaluated.
However, the main theorem in this section shows that norm of   $S_t m$ compares naturally with that of $m_+$ and $m_-$.

\begin{Thm}\label{Mult}
  Let $p$ be even and let $m :\R^n\to \C$ be a $(p,p)$ Fourier multiplier supported in  $\R^{n-1}\times [-A,A]$.
  Suppose that, for every $t>0$,  $T_{S_t m} \in \mathcal R$. Then,
for every  $t\ge t_0=\frac{(p-2)A}{4}$,
\begin{equation}\label{main}
  |||S_t m |||_{p,p}  \leq  {p\choose \frac p2}  ^{\frac 1p} \left( |||m_+|||_{p,p}  |||m_-|||_{p,p} \right)^{\frac 12}.  \end{equation}
%Equality holds when   $x'\to m(x', x_n)$ is even   for every   $x_n\in\R $  and
 % $x_n\to m(x', x_n)$ is even   for every $x'\in\R^{n-1}$, and    $p=2$.
\end{Thm}

\medskip
\noindent
{\it Remark.} If $m$ is real-valued,
it is not too difficult to prove that
$T_m\in {\cal R} $   if and only if $x'\to m(x', x_n)$  and
$x_n\to m(x', x_n)$ are both even.
 For such $m$'s,
 $|||  m_+|||_{p,p} = |||  m_-|||_{p,p} $ and
 \eqref{main} reduces to
 \BEL{mainR}|||S_t m |||_{p,p}^\R = |||S_t m |||_{p,p}  \leq  {p\choose \frac p2}  ^{\frac 1p} |||m_+|||_{p,p}. \EE Since $2 =\left(\sum_{k=0}^p {p\choose k}\right)^{\frac 1p} > {p\choose \frac p2} ^{\frac 1p}
$, \eqref{mainR} is an improvement over the trivial estimate
$||| S_t m |||_{p,p}   \leq 2|||\tau_{t} m_+|||_{p,p} = 2|||  m_+|||_{p,p}.$
Also, observe that   $|||m|||_{p,p} = |||m|||_{p',p'} $, and so
\eqref{main} is equivalent to
\begin{equation}\label{e-dualp}|||S_t m |||_{p',p'}  \leq  {p\choose \frac p2}  ^{\frac 1p} \left(|||m_+ |||_{p',p'} |||m_- |||_{p',p'} \right)^{\frac 12} .\end{equation}
\medskip
\noindent
{\it Proof.} Let $T_{S_t m} = T_t = T_+ + T_-$, where $T_+ = T_{\tau_t m_+}$ and $T_{-} = T_{\tau_{-t} m_-}$.
Since $T_{S_t m} \in {\cal R}$,  $|||S_t m |||_{p,p} =  |||S_t m |||_{p,p}^\R =\sup ||T_{t} f||_p$,
where $||f||_p =1$. By density, we can fix such an $f\in C_0^\infty(\R^n )$, and estimate $||T_{t} f||_p $.
 Then
\begin{eqnarray}\nonumber
 T_+ f(y)
&=& \int_{\R^{n }}     m_+ (x',\,x_n-t)  \hat f(x) e^{2\pi i xy}dx
\\\nonumber   &=&  \int_{\R^{n }}      m_+ (x)\widehat f(x',\,x_n+t) e^{2\pi i ((x_n+t)y_n + x'y')}dx
\\\label{T+}   &=
  &  e^{ 2\pi i ty_n}\int_{\R^{n }}    m_+ (x) {\cal F}(fe^{-2\pi ity_n })  (x) e^{2\pi i xy}dx
 \end{eqnarray}
where ${\cal F}$ is the Fourier transform on $\R^{n }$.
Define $\psi = \psi_+ +\psi_-$ by $\psi_+(x)= m_+(x){\cal F} (fe^{-2\pi ity_n }) (x)$,
 and $\psi_-(x)= m_-(x){\cal F} (fe^{2\pi ity_n })  (x)$. By \rq{T+}, $T_+ f(y)= e^{  2\pi i ty_n} \widehat {\psi_+}(-y)$. A similar calculation shows that \begin{align}\nonumber
T_- f(y) &=
e^{ -2\pi i ty_n}\int_{\R^{n }}    m_- (x){\cal F}(fe^{2\pi ity_n }) (x) e^{2\pi i xy}dx \\ \label{e-T-} &= e^{  -2\pi i ty_n} \widehat {\psi_-}(-y).
\end{align}

By \eqref{FofSplit}, $T_t f(y)=  e^{  2\pi i ty_n} \widehat {\psi_+}(-y) + e^{  -2\pi i ty_n} \widehat {\psi_-}(-y)= \widehat{S_t \psi}(-y)$.
By hypothesis $T_t f$, and hence $\widehat{S_t\psi}$, are real-valued.
Since $m$ is supported  in $\R^{n-1}\times [-A, \ A]$,
so is  $\psi(x)$.
Thus, we can apply Lemma \ref{CWt} with $F=\psi$ and
by \eqref{Newt}
and H\"older's inequality,
\begin{eqnarray}\nonumber|| \widehat {S_t \psi } ||_p &\equiv& {p\choose \frac p2}^{\frac 1p} \left( \int_{\R^n} \left(  \widehat  \psi_+ (y , t)  \right)^{\frac p2} \left(
\widehat  \psi_-  (y ,t)\right)^{\frac p2}dy\right)^{\frac 1p} \\ \label{EndMult} &\le& {p\choose \frac p2}^{\frac 1p} || \widehat  \psi_+ ||_p^{\frac 12}|| \widehat  \psi_- ||_p^{\frac 12}.
\end{eqnarray}
But $|| \widehat  \psi_+ ||_p =||T_+  f  ||_p \le |||\tau_t m_+ |||_{p,p}   =|||m_+ |||_{p,p} $
 and likewise $|| \widehat  \psi_- ||_p \le |||m_- |||_{p,p} $,
providing the required
estimate on $||T_{t} f||_p = || \widehat {S_t \psi } ||_p$, and proving \eqref{main}.
$\Box$

\medskip
Our next propositions
estimate $|||m_+|||_{p,p} $ in terms of $|||m|||_{p,p} $ and $c_p $.
Similar results hold for $|||m_-|||_{p,p} $ with similar proofs.
 \medskip
 \begin{Prop}\label{m+}
 Let $m:\R^n \to\C$ be a  $(p,p)$ Fourier multiplier.   Then,
 \begin{equation}\label{e-m+}
 |||m_+|||_{p,p}  \leq c_p   |||m|||_{p,p} .
 \end{equation}

If  $T_m\in \mathcal R$   we also have
\begin{equation}\label{e-rem+}
 |||m_+|||_{p,p}^\R \leq c_p ^\R  |||m|||_{p,p}^\R.
 \end{equation}

If $p=2$ and  $m$ is real-valued and even with respect to $x_n$, then equality holds in \rq{e-rem+}.

\end{Prop}
 \medskip
 \noindent{\it Proof.} Let $H(x)=\chi_{\{x_n>0\}}(x)$ be the half-plane multiplier so that $m_+ = mH$.
But $|||mH|||_{p,p} \leq  |||m|||_{p,p} |||H|||_{p,p}$ and \rq{e-m+} follows if
$|||H|||_{p,p} \le |||s|||_{p,p} = c_p$. To see this, note that
\begin{align*} T_H f(y)&=\int_{\R^n}  s(x_n) \hat f(x_n, x')  e^{2\pi i y x}dx  \\
&= \int_{\R } s(x_n) e^{2\pi i y_nx_n} {\cal F}_1 f( y', y_n) dx_n = T_{s} f(y',\cdot) .\end{align*} So, $\int_{\R}|T_H f|^p dy_n \le |||s|||_{p,p}^p ||f(y',\cdot)||_{L^p(\R)}^p$, and integration over $y'$
proves that $|||H|||_{p,p}\leq c_p$.
Equality can  be proved by setting $f=g(x') h(x_n)$, so that $T_H f(y)= g(y') T_sh(y_n)$.
The proof of \eqref{e-rem+} is similar.

To finish the proof,
we must show that $|||m|||_{2,2}^\R \le \sqrt 2|||m_+|||_{2,2}^\R$.
The standard proof of the identity $|||m|||_{2,2}=||m||_\infty$ can be easily modified to show that
$|||m|||_{2,2}^\R = ||m||_\infty = \sup  ||T_{m} f||_2$,   where the sup is taken over  the $f\in L^2(\R^n)$ that are even in $x_n$, with $||f||_2=1$. Fixing such an $f$,
and using the Plancherel theorem,
$$||T_m f||_2  = \sqrt{2} ||T_{m_+} f||_2  \leq   \sqrt{2} |||{m_+} |||_2^{\R}.
 \quad \Box
$$

 \medskip
Our next theorem provides a lower bound for $||S_t m||_{p,p}$ in terms of $c_p$.
Since $\dsize \lim_{p\to 1}c_p =\lim_{p\to\infty} c_p =\infty$, inequality \eqref{e-upperSt} shows that $|||S_tm|||_{p,p}$ cannot be  bounded above by a constant independent of $p$. Here, we assume $\dsize \ell_+  = \lim_{x\to 0^+} m(x)$ and $\dsize \ell_-  = \lim_{x\to 0^-} m(x)$ exist and that $\ell =\max\{|\ell_+|, |\ell_-| \} >0$.
\medskip
\begin{Thm}\label{HalfLine}
Let   $m\in L^\infty(\R,\,\C)$ be  a $(p,p)$ Fourier multiplier.
Then, for every $t>0$ and every $1<p<\infty$,  %
\begin{equation}\label{e-upperSt}
 |||S_t m |||_{p,p} \ge   \ell c_p.
 \end{equation}
\end{Thm}

\noindent
{\it Remark:} Suppose $m$ is supported in $[-A,A]$,  $T_{S_t m}\in \mathcal R$ for every $t>0$, $p$ is even and $t\ge t_0=\frac{(p-2)A}{4}$. Then
by Theorem \ref{HalfLine} and Theorem \ref{Mult},
\begin{equation}\label{St2way}
\ell\, c_p\leq   |||S_t m |||_{p,p}  \leq    c_p {p\choose \frac p2}^{\frac 1p}  |||m |||_{p,p}.
\end{equation}
\medskip
\noindent
{\it Proof.} Fix $m$, $t$ and $p$. Without loss of generality,   $\dsize \ell =\ell_+ = \lim_{x\to 0^+} m(x) >0$.
The  norm of a  multiplier is invariant by translation and dilation,
so the $(p,p)$ norm of $S_t m$ is the same as that of $k(x) = \tau_{-t}S_tm(x)= m_-(x+ 2 t)+ m_+(x)$. For $\lambda>0$, let $k_\lambda(x)= k(x/\lambda)$. Note that
$\dsize \lim_{\lambda\to\infty} k_\lambda(x)= \ell  s(x)$.
Fix $f\in C^\infty_0(\R,\C)$. By the Lebesgue dominated convergence theorem,
\begin{eqnarray*}
\lim_{\lambda\to\infty} T_{k_\lambda} f(y) &=& \lim_{\lambda\to\infty} \int_{\R} \hat f(x) k_\lambda(x) e^{2\pi i y x} dx \\ &=&  \ell \int_{\R } \hat f(x) s(x) e^{2\pi i y x} dx    =  \ell  T_{s} f(y). \end{eqnarray*}
By Fatou's Lemma
\begin{eqnarray*} \liminf_{\lambda \to +\infty} \int_\R |T_{k_\lambda} f(y)|^pdy  &\ge&
\int_\R \liminf_{\lambda \to +\infty}|T_{k_\lambda} f(y) |^pdy  \\ &=& \ell^p \int_\R |T_{s} f(y)|^pdy  =  \ell^p\,||T_{s } f||_p^p.
\end{eqnarray*}
We obtain
\BEL{11}
\ell c_p = \ell \,|||s|||_{p,p}  \leq  |||k_\lambda|||_{p,p} = |||S_{t }m|||_{p,p} .
\ \ \Box
\EE

 \medskip
The proof of Theorem  \ref{HalfLine} also provides  a lower bound (and in some cases an exact value)  for
 the $(p,p)$ norm of $m_+$.
\begin{Cor}\label{P-extram+}    Let  $m:\R\to\C$ be a  $(p,p)$ multiplier, continuous at $x=0$, with $\ell = m(0) \ne 0$. Then,
\begin{equation}\label{e-2cmplxm+}
 c_p|\ell|  \leq |||m_+|||_{p,p} \leq c_p|||m|||_{p,p}.
  \end{equation}
Suppose also that $m$ is real-valued and even. Then
  \begin{equation}\label{e-2realm+}
 c_p^\R  |\ell |  \leq |||m_+|||_{p,p}^\R  \leq c_p^\R  |||m|||_{p,p}.
  \end{equation}
The four inequalities in \rq{e-2cmplxm+} and \rq{e-2realm+} are equalities if $\hat m \in L^1(\R)$ and $\hat m \ge 0$.
 \end{Cor}

\noindent
{\it Remarks.} The same bounds hold for $|||m_-|||_{p,p}$ and $|||m_-|||_{p,p}^\R$.
Equality also holds in \eqref{e-2cmplxm+} if  $|\ell| =||m||_\infty $ and $p=2$ because $||m||_\infty = |||m|||_{2,2}$.
It also holds in \eqref{e-2realm+}, since $T_m\in \mathcal R$ implies $|||m|||_{2,2}^\R =|||m|||_{2,2}$.
  \medskip
 \noindent
   {\it Proof of Corollary \ref{P-extram+}.} The upper bounds on the norms of $m_+$ are in Proposition \ref{m+}. The proof of \eqref{e-upperSt} establishes   the lower bounds.  To prove that equality holds in \eqref{e-2cmplxm+},
   assume that $\hat m \in L^1(\R)$ and $\hat m \ge 0$  so that  $||\hat m||_1 = \int_\R \hat m(x)dx = m(0)=\ell$.
   For every $f\in C^\infty_0(\R^n)$,  $T_m f(x)= f* \hat m(-x)$, so by Young's inequality,
   $$||T_m f||_p\leq ||f||_p ||\hat m||_1= ||f||_p \ell.$$ Thus $|||m|||_{p,p}\leq \ell$, which implies equality in \eqref{e-2cmplxm+};  the proof for \eqref{e-2realm+}  is similar. $\Box$

   \medskip
If $m $ is increasing on $\R$ with $\sup m(x)=1$ and $\inf m(x)=0$, then $|||m|||_{p,p} = c_p$ (see \cite{DL2}).
Such a multiplier cannot have  compact support, but the one in   the example below  does, again with $|||m_+|||_{p,p}=c_p$.

 \medskip
 \noindent
 {\it Example}. Let $m(x)= (1-|x|)\chi_{(-1, 1)}(x)$. It is easy to verify that  $\hat m (y)= \frac{\sin ^2(\pi  y)}{  \pi ^2 y^2} \ge 0$, and
  $||m||_\infty=\ell = 1$.  Corollary \ref{P-extram+} shows  that $|||m_+|||_{p,p}= c_p$ and $|||m_+|||_{p,p}^\R= c_p^\R$, and $|||m|||_{p,p} =1$.

\medskip
 As a geometric application, we estimate the multiplier norms of "split polygons".
 Let $P$ be a polygon with diameter $2A$ in $\R^2$, and assume that $\chi_{P }(x_1, x_2)$ even in both $x_1$ and $x_2$.
It is well known that $\chi_P$ is a $(p,p)$ multiplier for every $1<p<\infty$, with bounds on its multiplier norm that depend on $p$ and the number of sides (see \cite{Co1}). Let $P_+=   P  \cap\{x_2>0\}$. Assume that $p$ is even, and    $n_p$ is defined by \eqref{segment}.

\begin{Prop}\label{T-Poly} If the intersection of $P$ and the line $\{x_2=0\}$
is a segment $I$,
and $t\ge t_0=\frac{(p-2)A}{4}$, then
\begin{equation}\label{e-poly}
n_p c_p \leq ||| S_t (\chi_{P})|||_{p,p} \leq {p\choose \frac p2}^{\frac 1p} |||\chi_{P_+}|||_{p,p}.
\end{equation}
\end{Prop}

\medskip
\noindent
{\it Proof.}
Since $P$ is symmetric, $S_t (\chi_P)\in {\cal R}$ and the second inequality follows from  Theorem \ref{Mult} (see also the Remark following it).  For the first, we can assume by dilation that $I = (-1,1)$.
The  $(p,p)$ norm of  $S_t (\chi_{P})$ is the same as that of $ \tau_{-t}S_t (\chi_{P})(x)$, where the translation acts only on the $x_2$ variable.
Let $k(x) =\tau_{-t}S_t(\chi_{P})(x)$. For $\lambda>0$, let $k_\lambda(x)= k(x_1,\,\lambda^{-1}x_2 )$.
Note that $ k_\lambda(x)\equiv 0$  if $-2t\lambda< x_2 < 0$,
and $k_\lambda( x_1,\ 0) \equiv 1$ whenever $-1<x_1<1$, and is $\equiv 0$ otherwise.  Thus
 \BEL{limit} \dsize \lim_{\lambda\to\infty} k_\lambda(x_1, x_2)=   s(x_2)\cdot \sigma(x_1)
 \EE
where $\sigma(x_1)= \chi_{(-1,1)}(x_1)$. By Lebesgue's Theorem and Fatou's Lemma (as in the proof of Theorem \ref{HalfLine}), we get
 $$n_p c_p =  |||s\cdot \sigma|||_{p,p}  \leq \liminf_{\lambda \to +\infty} |||k_\lambda |||_{p,p} = |||k |||_{p,p} = |||S_{t}(\chi_P)|||_{p,p}.  \ \ \Box$$

\medskip
\noindent
{\it Example.} Let $Q$ be the square with corners at the points $(\pm A,0)$  and $(0, \pm A)$.
Assume $p$ is even and $t\ge t_0=\frac {(p-2)A }{4}$.
Since $|||\chi_{Q_+} |||_{p,p}\le c_p^3$ (see \cite{De}), Proposition \ref{T-Poly} implies
\BEL{2tri}
n_p c_p\leq  |||S_t (\chi_Q) |||_{p,p} \leq   {p\choose \frac p2}^{\frac 1p}    c_p^3.
\EE

 \end{document}